\documentclass[12pt,a4paper]{article}
\usepackage[english]{babel}
\usepackage{amssymb}
\usepackage[dvips]{graphicx}

\def\C{{\mathbb C}}

\def\f{\varphi}

\def\l{\lambda}

\def\s{\sigma}
\def\Z{{\bf Z}}

\def\w{{\omega}}
\def\ww{{\overline \w}}
\def\LL{{\mathcal L}}

\def\F{{\mathcal F}}

\def\R{{\mathcal R}}
\def\U{{\mathcal U}}

\def\beq{\begin{equation}}
\def\eeq{\end{equation}}

\def\a{\alpha}
\def\b{\beta}
\def\g{\gamma}
\def\G{\Gamma}

\def\k{\kappa}
\def\we{\wedge}
\def\wt{\widetilde}

\def\e{\varepsilon}
\def\ov{\overline}

\def\ZZ{{\mathbb Z}}

\renewcommand{\l}{\lambda}

\newtheorem{thm}{Theorem}[section]
\newtheorem{pr}[thm]{Proposition}
\newtheorem{cor}[thm]{Corollary}

\newtheorem{dfn}[thm]{Definition}

\begin{document}

\title{{\bf Automorphisms of $P_8$ singularities \\ and the complex crystallographic groups}}

\author{Victor Goryunov and Dmitry Kerner}
\date{}
\maketitle

\begin{abstract}
The paper completes the study of symmetries of parabolic function singularities
with relation to complex crystallographic groups that was started in \cite{GM,X9}. 
We classify smoothable automorphisms of $P_8$ singularities 
which split the kernel of the intersection form on the second homology. For such automorphisms,
the monodromy groups acting on the duals to the eigenspaces with degenerate intersection form are then identified as some of 
complex affine reflection groups tabled in \cite{P}.
\end{abstract}

Singularity theory has been maintaining close relations with reflection groups since its early days, starting with the famous
works of Arnold and Brieskorn \cite{A,Br}. The first to emerge, as simple function singularities, were the 
Weyl groups $A_k, D_k$ and $E_k$. They were followed by the $B_k$,  $C_k$ and $F_4$ as simple boundary 
singularities or, equivalently, functions invariant under the involution \cite{Ab,Slo}. Since then appearance of Weyl groups 
in a singularity problem became a kind of criterion of naturalness of the problem \cite{AGV2,AGLV1,AGLV2}. 

The next to receive a singularity realisation was the list of all Coxeter groups, in the classification of 
stable Lagrangian maps \cite{Gi}. Some of the
Shephard-Todd groups appeared in \cite{G1,G2,GB} in the context of simple functions equivariant with respect to a cyclic
group action. And finally, the first examples of complex crystallographic groups came out in \cite{GM,X9} 
in connection with the symmetries of parabolic functions $J_{10}$ and $X_9$. 
The affine reflection groups appeared there as monodromy groups, which 
 is similar to the first realisations of other classes of reflection groups. This time it was the equivariant monodromy
corresponding to the symmetry eigenspaces $H_\chi$ in the vanishing second homology on which the intersection form 
$\s$ has corank 1.

This paper completes the study of cyclically equivariant parabolic functions started in \cite{GM,X9} and considers the $P_8$ singularities.
The approach we are developing here is considerably shorter, with less calculations and more self-contained. This is allowed by 
a preliminary observation that the modulus parameter may take on only exceptional values if corank$(\s|_{H_{\chi}}) = 1$: the
$j$-invariant of the underlying elliptic curves must be either 0 or 1728 (see Section \ref{Sker}).

\medskip
The structure of the paper is as follows.

Section 1 describes the crystallographic groups which will be involved.

Section 2 gives a classification of smoothable cyclic symmetries of $P_8$ singularities possessing eigenspaces $H_{\chi}$
with the property as already mentioned.

In Section 3, we obtain distinguished sets of generators in such eigenspaces and calculate the intersection numbers of the
elements in the sets. This allows us to describe the Picard-Lefschetz operators generating, as complex reflections,  the equivariant
monodromy action on the $H_{\chi}$.

Finally, in Section 4, we pass to the vanishing cohomology. In its subspaces dual to the $H_{\chi}$, we consider hyperplanes
formed by all the cocycles taking the same non-zero value on a fixed generator of the kernel of the hermitian intersection form
$\s|_{H_{\chi}}$. We show that the equivariant monodromy group acting on each of these hyperplanes is a complex 
crystallographic group.
 
\medskip
Both authors gratefully acknowledge the support of the EC Marie Curie Training Site (LIMITS) on Algebraic Geometry, 
Singularities and Knot Theory which funded the visits of the second author to Liverpool.
The research of the second author was also supported by the Skirball postdoctoral fellowship of the Center 
of Advanced Studies in Mathematics, Be'er Sheva, Israel.


\section{The crystallographic groups}
\begin{figure}[hbtp]
\begin{center}
\includegraphics[width=0.7\textwidth]{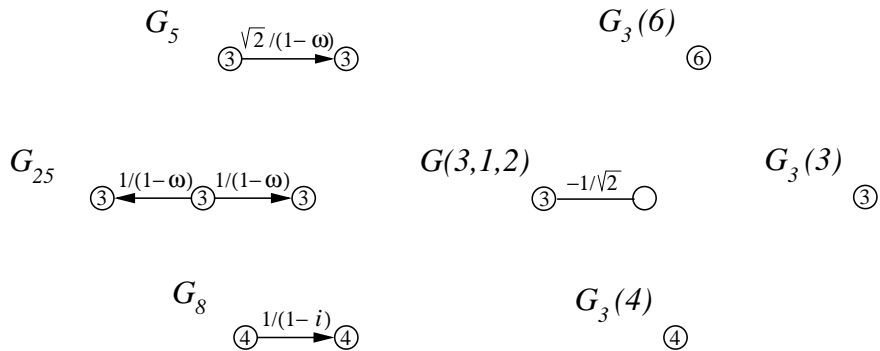}
\caption{Linear parts of the crystallographic groups. Each vertex represents a unit root. The order of the corresponding reflection 
is written inside the vertex (order 2 omitted). An egde $a \to b$ is equipped with the hermitian product $\left< a, b \right>$ 
(no orientation is needed when the product is real). There are no edges between orthogonal roots.}
\label{Fgr}
\end{center}
\end{figure}

We remind the description of the complex affine reflection groups 
which will be relevant to our singularity constructions. In the notations of \cite{P}, the groups
are
$$
[K_3(3)], [K_3(4)], [K_3(6)], [G(3,1,2)]_2,  [K_5], [K_8], [K_{25}]\,.
$$ 
The linear parts of these groups are Shephard-Todd groups $L$  shown in Figure \ref{Fgr}, in the related notation of \cite{ST}. 
There and in what follows $\w=exp(2\pi i/3)$.
The rank of $L$ is the number of the vertices in its diagram. 
Each of the groups involved is a semi-direct product of its linear part and translation lattice. 
The lattice of $[G(3,1,2)]_2$ is spanned over $\ZZ$ by the
$G(3,1,2)$-orbit of any order 2 root, while for any other group it is spanned the $L$-orbit of any root. For example, our list contains
all three one-dimensional groups, the $[K_3(m)]$, each generated by an order $m$ rotation in $\C$ and the lattice $\ZZ[1,i]$ if $m=4$ or 
$\ZZ[1,\w]$ if $m=3,6$.

\section{The symmetries}

\subsection{Cubic curves}
Our first step in the study of linear automorphisms of the surfaces
$$x^3+y^3+z^3+3\a xyz =0\,, \qquad \a^3 \ne -1,$$
in $\C^3$, is consideration of the projective version of the question. This is a
description of automorphisms of projective curves $C_\a$ given by the same
equations in $\C P^2$. The question is classical, but, for the sake of the exposition,
we indicate a possible elementary approach.

A curve $C_\a$ has nine inflection points. There are four triplets of lines
in $\C P^2$ passing through these points, each line containing exactly three inflections
and each triplet containing all nine inflections. The lines are zero sets of the linear forms
$$
\begin{array}{lll}
x, & y, & z; \\
x+y+z, \quad & x+\w y + \ww z, \quad & x + \ww y + \w z; \\
x+y+\w z, \quad & x + \w y + z, \quad & \w x + y + z; \\
x+y+\ww z, \quad & x + \ww y + z, \quad & \ww x + y + z\,.
\end{array}
$$
A projective automorphism of $C_\a$ permutes the triplets and reorders the lines
within the triplets. It is sufficient to consider the images of the first triplet
only. A routine study of all the possibilities, followed by diagonalisation of each 
projective transformation and rewriting the equation of an appropriate curve in the 
projective eigencoordinates, give us 

\begin{pr}\label{Pproj}
Any projective automorphism of a cubic curve can be reduced, by a choice of projective
coordinates, either to the diagonal projective transformation 
$(x:y:z) \mapsto (\k_xx:\k_yy:\k_zz)$ 
of one of the following
cubics $f=0$ or to the inverse of such a transformation:
$$
\begin{array}{|l|l|l|l|} 
\hline 
&&& \vspace{-8pt} \\
 & \qquad\qquad\qquad  f & \k_x:\k_y:\k_z & 
{\rm monomial \  basis \  of \  the \  local \  ring \  of \  } f
\\ 
&&& \vspace{-8pt} \\
 \hline &&& \vspace{-8pt} \\
1 & x^3 + y^3 + yz^2   & \,\, \w : 1 : -1 & 
1, x, y, z, xy, xz, y^2\sim z^2, xy^2 \sim xz^2 \\ &&& \vspace{-8pt} \\
2 & x^3 + y^3 + yz^2   \simeq x^3 + y^3 + z^3 & \,\, \w : 1 : 1 &
1, x, y, z, xy, xz, y^2\sim z^2, xy^2 \simeq xyz \\ &&& \vspace{-8pt} \\
3 & x^2y + y^2z + z^2x & \,\, 1 : \w : \ww &
1,x,y,z,x^2 \sim yz, y^2 \sim xz, z^2 \sim xy, xyz  
\\ &&&  \vspace{-8pt} \\
4 & x^2z + xy^2 + z^3 & \,\, 1 : i : -1 &
1,x,y,z,yz,x^2\sim z^2,xz\sim y^2,x^3 \\  &&& \vspace{-8pt} \\
5 & x^3 + y^3 + yz^2 + 3 \a xy^2\,, \,\, \a^3 \ne -\frac14 & \,\, 1:1:-1 & 
1, x, y, z, xy, xz, y^2\sim z^2, xy^2 \sim xz^2 \\ &&& \vspace{-8pt} \\
6 & x^3+y^3+z^3+ 3 \a xyz\,, \,\, \a^3 \ne -1  & \,\, 1:\w:\ww  &
1,x,y,z,xy,xz,yz,xyz \\ &&& \vspace{-8pt} \\
7 & x^3+y^3+z^3+3\a xyz\,, \,\, \a^3 \ne -1 & \,\, 1:1:1 &
1,x,y,z,xy,xz,yz,xyz \vspace{-8pt} \\ &&& \\ \hline
\end{array}
$$
\end{pr}

\medskip
The contents of the last column of the table will be used in Section \ref{Ssmosym}.
 
\subsection{Kernel character constraints}\label{Sker} 
We now lift a projective automorphism of a cubic curve $f=0$ to a linear
transformation $g$ of $\C^3$. It multiplies $f$ by a non-zero constant 
and hence acts on the target $\C$ of the function $f$. 
We consider deformations of the function-germ $f: (\C^3,0) \to (\C,0)$ 
as usual: choosing a small ball $B$ in $\C^3$ centred at the origin and deforming $f$ within the ball. However, 
we restrict our attention to deformations which are $g$-equivariant with respect to the actions of $g$ on
$\C^3$ and $\C$. 
 
\begin{dfn} {\em
If the function-germ $f$ admits $g$-equivariant deformations 
with non-singular zero levels, we say that the symmetry $g$ of $f$
is} smoothable.
\end{dfn}

In this case, the smooth zero level $V \subset B$ of a small $g$-equivariant 
deformation of $f$ is a Milnor fibre of $f$, and $g$ acts on $V$ and on its 
second homology. Assuming from now on the order $r$ of $g$ finite, we obtain the 
splitting into a direct sum of character subspaces 
\begin{equation}
H_2(V,\C) = \oplus_{\chi^r=1} H_{2,\chi}
\label{Esum}
\end{equation}
so that $g$ acts on each summand as multiplication by the relevant 
$r$th root $\chi$ of unity.

We extend the intersection form $\sigma$ from the lattice $H_2(V,\ZZ)$ to
the hermitian form $\sigma_h$ on $H_2(V,\C)$. The symmetric intersection form 
of a parabolic function singularity in 3 variables is negative semi-definite, with a two-dimensional kernel. 

\begin{dfn} {\em
We say that a smoothable symmetry $g$} splits {\em the kernel of
$\sigma_h$ if this kernel is shared by two different character subspaces
$H_{2,\chi}$. If the restriction of $\sigma_h$ to a $H_{2,\chi}$ is 
degenerate, we call $\chi$ a} kernel character.
\end{dfn}     

The aim of our construction, a complex crystallographic group, will be acting 
on a hyperplane $\G$ in the cohomology subspace $H_{\chi}^2 = \{\gamma \in H^2(V,\C): g^*(\gamma) = \chi\gamma\}$ where
$\chi$ is a kernel character (cf. \cite{GM,X9}).
For this, our symmetry $g$ should split the kernel of $\sigma_h$, and the hyperplane $\G$ is defined as 
the set of all elements of $H_{\chi}^2$ taking the same non-zero value
on a fixed generator of the line $Ker(\sigma_h) \cap H_{2,\chi}$.

\medskip
Our next step is to show that there is no kernel splitting in any of the modular
cases of Proposition \ref{Pproj}. 

Assume $g$ is a lifting of one of the projective automorphisms from the table there.
Take a monomial basis $\f_0=1, \f_1, \dots, \f_7$ of the local ring
of $f$, with $\f_7$ of degree 3.
Consider the sub-unfolding of the $\R$-miniversal unfolding of $f$ in
which the modulus parameter does not participate:
$$
\F: \C^3 \times \C^{6} \to \C \times \C^{6}\,, 
\quad 
(x,y,z,\l_1,\dots,\l_{6}) \mapsto
(f+\l_1\f_1+\dots\l_{6}\f_{6},\l_1,\dots,\l_{6})\,.
$$ 
This map is defined globally. Its smooth fibres are Milnor fibres of $f$.
The action of $g$ naturally extends to the source and target of the map,
making $\F$ $g$-equivariant. We will denote the target of $\F$ by $\U$,
points there will be $u=(u_0,\dots,u_{6})$.

We simultaneously compactify all fibres $V_u$, 
$u \in \U$, of $\F$ adding to each of them a copy 
$C_u = \ov V_u \setminus V_u$ of the same projective curve $f=0$. 
Each $\ov V_u \subset \C P^3$ is smooth near $C_u$. 

For non-critical values of $u$, take the intersection homomorphism
$\s:  H_2(V_u;\ZZ) \to H^2(V_u;\ZZ)$.
The Leray map sends $H_1(C_u; \ZZ)$ isomorphically onto the 2-dimensional 
kernel of $\s$ \cite{L}.

Consider the 2-form $w = dx \we dy \we dz \we d\l_1 \we \dots \we d\l_{6}/ 
du_0 \we \dots \we du_{6}$ on the fibres
$V_u$. According to \cite{L}, its residue $\b$ along $C_u$ is a non-zero
holomorphic 1-form which does not depend on $u$.

The action of the automorphism $g$ on a $g$-symmetric fibre $V_u$ 
extends to $\ov V_u$.  The residue map
$Res: H^2(V_u;\C) \to H^1(C_u;\C)$ is $g$-equivariant: if $g^*$ multiplies a 2-form
by $\chi$, then it multiplies the residue of the form by $\chi$ too. Therefore, if $g$
splits $Ker(\s_h)$ on $V_u$, then its action on $H^1(C_u;\C)$ must also
have two distinct conjugate eigenvalues, with $\b$ and $\ov \b$ as eigenvectors. 
Hence the curve $C_u=\{f=0\}$ must be either $\C/\ZZ[1,i]$ or $\C/\ZZ[1,\w]$. 
Respectively, the kernel characters are either $\pm i$ or 
$(\w,\ov \w)$ or $(-\w,-\ov\w)$.  

Thus we have:

\begin{pr}
A smoothable automorphism $g$ splits the kernel of the intersection form
if and only if the section $w=dx \we dy \we dz /df$ of the vanishing cohomology
fibration of $f$ is a $g^*$-eigenvector with the eigenvalue from the set
$\{\pm i, \pm \w, \pm\ov\w\}$. The eigenvalue and its conjugate are the 
kernel characters. 
\end{pr}

If $g^*(w) = \chi w$ then for the hyperplane $\G$
mentioned above  one can take $Res^{-1}(Res(w)) \cap H_{\chi}^2$.

Since the degree of $f$ is 3, the eigenvalue of $g^*$ on $w=dx \we dy \we dz /df$
depends only on the projectivisation of $g$. 

\begin{cor}\label{Clift}
All smoothable liftings of the first four automorphisms from Proposition
\ref{Pproj} split $Ker(\s_h)$. On the other hand, none of the liftings of the
last three cases does the same. 
\end{cor}

\noindent
{\bf Remark.} The argument can be easily modified to the 
symmetries of $X_9$ and $J_{10}$ function singularities, thus explaining
the absence of the moduli in the classification tables in \cite{GM,X9}.

\subsection{Smoothable symmetries}\label{Ssmosym}
Assume we have two actions of a finite cyclic group, on $(\C^k,0)$ and on $(\C,0)$.
Consider two function-germs $f_1,f_2:  (\C^k,0) \to (\C,0)$ equivariant with
respect to these actions: $f_i \circ \rho = \rho \circ f_i$ where $\rho$ is a generator of
the group. We say that the two functions are {\em $\R_\rho$-equivalent\/} if one
can be transformed into the other by a $\rho$-equivariant diffeomorphism-germ 
of $(\C^k,0)$. 

We apply this notion in the context of our cubic function $f$, its {\em diagonal\/}
symmetry $g$ and the induced action of $g$ on $\C$.  With a minor abuse of the
language, we will still be calling the corresponding equivalence the $\R_g$-equivalence.
For example, an $\R_g$-miniversal unfolding of $f$ 
is the restriction $F$ of the unfolding $\F$ of the previous subsection to 
$\F^{-1}(U)$, where $U \subset \U$ is the set of fixed points of the natural
action of $g$ on $\U$. 
Equivalently, in this case, for an $\R_g$-miniversal deformation of 
$f$ one can take $f +$ arbitrary linear combinations of the elements $\psi_1, \dots, \psi_\tau$ of 
a monomial basis of the local ring of $f$ that are multiplied by $g$ by the
same factor as $f$. The number $\tau$ appearing here will be called the 
{\em $\R_g$-codimension\/} of $f$.

\begin{pr}\label{Proff0}
Assume a symmetry $g$ is smoothable. Then $g$ multiplies $f$ by the same
factor as it multiplies one of the monomials $1, x, y, z$.
\end{pr} 

Indeed the conclusion is equivalent to generic fibres of an $\R_g$-miniversal unfolding
not being singular at $0 \in \C^3$. 
 
With the help of Corollary \ref{Clift} and the last Proposition, we obtain after
straightforward calculations which we prefer to omit:

\begin{thm}\label{Tclass} 
The complete list of  smoothable symmetries of $P_8$ function singularities 
on $\C^3$ which split the kernel of the intersection form is given in Table \ref{Ta}. 
Our classification is up to a choice of generators of the same cyclic group.
\end{thm}

\begin{table}[p]
\begin{center}
\caption{Smoothable symmetries of $P_8$ singularities splitting $Ker(\s)$}
\vspace{-15pt}
$$\begin{array}{|c|c|c|c|c|c|c|}
\hline 
&&&&&&\\
f & g: x,y,z \mapsto & |g| & 
\begin{array}{c}
{\rm versal} \\
{\rm monomials}
\end{array}
&
\begin{array}{c}
\!\!\!{\rm kernel}\!\!\! \\
\chi
\end{array}
&
\begin{array}{c}
{\rm affine}\\
{\rm group}
\end{array}
 & 
{\rm notation} \\
&&&&&&\\
\hline 
\hline
&&&&&&\\
x^3    &   \w x,y,-z   &   6   &   1, y, y^2    &   -\w, -\ww  &  \!\! [K_5] \!\!  &   C_3^{(3,3)} \\
+y^3 &&&&&& \\
+yz^2 &   x, \ww y, -\ww z    &   6   &      1,x        &   -\w, -\ww   &  [K_3(6)]    &   (P_8|\ZZ_6)'   \\
&&&&&& \\
& \ww x, \w y, -\w z  &  6  & 1,xy  &  -\w, -\ww & [K_3(6)] & (P_8|\ZZ_{6})'' \\
&&&&&& \\
& -x, -\ww y, \ww z &  6  &  x  &  -\w, -\ww  &  -  &  (P_8/\ZZ_6)' \\
&&&&&& \\
& i\w x, iy, -iz  &  12  &  z  &  -\w, -\ww  &  -  &  P_8/\ZZ_{12} \\
&&&&&& \\
&  \w x, y, z  &  3  &  1,y,y^2,z  & \w, \ww & [K_{25}] & D_4^{(3)} \\
&&&&&& \\
&  x, \ww y, \ww z  &  3  &  1,x  & \w,\ww  &  [K_3(6)] & (P_8|\ZZ_{3})' \\
&&&&&& \\
&  \ww x, \w y, \w z  &  3   &  1,xy,xz & \w, \ww & [G(3,1,2)]_2 & P_8|\ZZ_3 \\
&&&&&& \\
\hline
&&&&&& \\
x^3 &  -x,-\ww y, -\ww z  &  6  &  x  & \w, \ww & - & (P_8/\ZZ_6)'' \\
+y^3 &&&&&& \\
+z^3 & -\w x, -y, -z  &  6  &  y,z  &  \w, \ww  &  [K_3(3)]  & P_8/\ZZ_6 \\
&&&&&& \\
\hline
&&&&&& \\
\begin{array}{c}
x^2y \\ + y^2z \\ + z^2x
\end{array}  
&  \ov\e_9 x, \w\ov\e_9 y, \ov{\w\e}_9 z  &  9   &  1 & \w, \ww & - & P_8|\ZZ_9 \\
&&&&&& \\
\hline
&&&&&& \\
x^2z   &  -x, -iy, z  &  4   &  1,z,z^2 & \pm i & [K_8] & C_3^{(2,4)} \\
+xy^2 &&&&&& \\
+z^3 & -\w x, -i\w y, \w z & 12 &  1  & \pm i  &  -  &  P_8|\ZZ_{12} \\
&&&&&& \\
& ix, -y, -iz  &  4  & x, yz  & \pm i & [K_3(4)]  &  P_8/\ZZ_4 \\
&&&&&& \\
& \ov\e_8 x , \e_8 y, -\ov\e_8 z  &  8  &  y  &  \pm i  &  -  &  P_8/\ZZ_8 \\
&&&&&& \\
\hline
\end{array}
$$
\label{Ta}
\end{center}
\end{table}

In the Table: the versal monomials are the monomials participating in $\R_g$-miniversal deformations 
(these are selected from the last column of the
table in Proposition \ref{Pproj} as those multiplied by $g^*$ by the same factor as $f$), 
the affine groups are those which will come out later as monodromy
groups acting on the hyperplanes in the character subspaces in the cohomology,
and $\e_k=exp(2\pi i/k)$.

\medskip\noindent
{\bf Remark.} During the proof of the Theorem, one obtains a few non-smoothable symmetries allowing linear terms in deformations:
$$
\begin{array}{|c|c|c|c|}
\hline 
&&&   \vspace{-8pt} \\
f & g: x,y,z \mapsto & |g| & 
{\rm versal \   monomials}    \vspace{-8pt}  \\
&&&\\
\hline 
\hline
&&&   \vspace{-8pt} \\
x^3 + y^3 + yz^2 & -\w x, -y, z  &  6  &  y   \vspace{-8pt} \\
&&&\\
x^2y + y^2z + z^2x  &  \w x, \ww y, z  &  3  &  x, xz    \vspace{-8pt}  \\
&&&\\
&  -\w x, -\ww y, -z  &  6  &  x    \vspace{-8pt}  \\
&&&\\
x^2z + xy^2 + z^3  &  x, iy, -z & 4  &  z, xz    \vspace{-8pt}  \\
&&& \\
& -ix, y, iz  &  4  &  x   \vspace{-8pt} \\
&&& \\
\hline
\end{array}
$$
It is possible to relate complex crystallographic groups to the two singularities
with 2-parameter miniversal deformations, but this will be done in another paper.   

\bigskip
Of course, one may consider a bit more general problem of finding all possible 
smoothable symmetry {\em groups\/} $G$ of parabolic singularities, which are not necessarily cyclic.
In this situation the second homology splits into irreducible represenations of $G$, and we may be still looking
for cases when the kernel of the intersection form is shared by two of them.
For $P_8$, for example, this means that $G$ should contain one of the symmetries $g$ of the Table, 
and hence the affine reflection group related to $G$ will be a lower rank subgroup of the one we are relating to $g$.
This does not leave too much dimensional room for further interesting crystallographic groups.


\section{Dynkin diagrams}
The equivariant monodromy of a $g$-equivariant function singularity $f$, that is, the one
within an $\R_g$-miniversal deformation of $f$ preserves the direct sum decomposition (\ref{Esum}).
Its action on an individual summand $H_{2,\chi}$ is generated by the Picard-Lefschetz operators $h_e$
corresponding to vanishing $\chi$-cycles $e$:
$$
h_e(c) = c - (1-\l_e) \left< c,e \right> e /\left< e, e \right>\,.
$$
Here $\l_e \ne 1$ is the eigenvalue of $h_e$: $h_e = \l_e e$. 

In this Section we obtain all the necessary information describing the monodromy on
the kernel character subspaces of the singularities of Table \ref{Ta}. We choose 
distinguished sets of generators of the subspaces $H_{2,\chi}$ (cf. \cite{Brmon,GZ, AGV2,AGLV1}), calculate 
the intersection numbers of the elements
of the sets, and find the eigenvalues $\l_e$. The data will be collected into Dynkin diagrams of the singularities.


\begin{pr}\label{Pdim}
If $\chi$ is a kernel character of a symmetric singularity from Table \ref{Ta},
then the rank of $H^2_{\chi}$, equivalently the rank of $H_{2,\chi}$, coincides
with the $\R_g$-codimension $\tau$ of the singularity.
\end{pr}

{\em Proof.} If $f$ is $g$-invariant, then a basis of one of the two $H^2_{\chi}$
is formed, in the notations of Section \ref{Sker}, by the sections
$\f_i w$, where the $\f_i$ are all $g$-invariant monomials within 
$\{\f_0, \dots, \f_7\}$. Their number is exactly $\tau$. (Cf. \cite{OS,W}.)

If $f$ is $g$-equivariant rather than invariant, the claim can be verified by
a direct case-by-case computation of the eigenvalues of $g^*$ on the 
sections $\f_0 w, \dots, \f_7 w$. The observation needs a case-free understanding.
\hfill{$\Box$}      

\medskip
Thus, if $\tau=1$ then the subspaces $H_{2,\chi}$ in the case of kernel $\chi$
are one-dimensional, hence contained in $Ker(\s_h)$, and therefore the monodromy
we are interested in is trivial. So, from now on we are considering only the cases
of at least two-parameter $\R_g$-miniversal deformations.

We will also forget about the $(P_8|\ZZ_3)'$ singularity. Indeed its symmetry $g$
is the inverse of the square of the $(P_8|\ZZ_6)'$ symmetry and the two 
equivariant miniversal deformations coincide. 
Hence the character subspaces with degenerate intersection form are
the same (only the actual character assignments differ: $\chi_3 = \chi_6^{-2}$)
and the monodromy is the same.  

\bigskip
For each remaining symmetry, we will choose a distinguished set of generators for $H_{2,\chi}$ 
in the way it has been done in \cite{G1,G2,GB}.
Namely, let $u_*$ be a point  in the complement $\C^{\tau} \setminus \Delta$ to the discriminant in the
base of the $\R_g$-miniversal deformation of our function. Set $V_{u_*}/g$ to be the quotient of the 
fibre $V_{u_*}$ under the cyclic
$g$-action and $(V_{u_*}/g)'$ its part corresponding to irregular orbits.  Denote
by $\pi$ the factorisation map $V_{u_*} \to V_{u_*}/g$. 
The homology group $H_2(V_{u_*}/g,(V_{u_*}/g)';\ZZ)$ is spanned by a 
distinguished set of relative vanishing cycles (such a set is defined in the traditional singularity theory manner). Denote them $c_1, \dots, c_k$.
The inverse image $\pi^{-1}(c_j)$ is the cyclic orbit of one of its components, 
let it be $\wt c_j$. For a character $\chi$, the chain
$$e_j=\sum_{i=0}^{order(g)-1} \chi^{-i}g^i(\wt c_j)$$
is an element of $H_{2,\chi}$.  The {\em vanishing $\chi$-cycles\/} 
$e_1, \dots, e_k$ span $H_{2,\chi}=H_{2,\chi}(V_{u_*},\C)$, not necessarily as a basis.
They are defined up to multiplication by powers of $\chi$ and the sign change.

\medskip
\noindent
{\bf Examples.}  (a) The function $x^n+y^2+z^2$ invariant under the transformation $g: (x,y,z) \mapsto (\e_n x,y,z)$ was denoted $A_1^{(n)}$
in \cite{G1}. 
The factorisation by the group action gives the boundary singularity $A_1: w + y^2 + z^2$, $w=x^n$, whose relative
vanishing homology is spanned by 
one semi-cycle $\{ w + y^2 + z^2 = 1: w,y,z \in {\mathbf R}, w \ge 0\}$. The corresponding vanishing $\chi$-cycle, $\chi^n=1$, $\chi \ne 1$, 
is formed by $n$ hemi-spheres with appropriate coefficients and has self-intersection $-n$ \cite{G1}, which is consistent
with the standard Morse vanishing cycle having self-intersection $-2$. 
The Picard-Lefschetz operator $h$ is the classical monodromy of the ordinary $A_{n-1}$ singularity. Its eigenvalue
on $H_{2,\chi}$ is $\l=\chi$, since the quasi-homogeneous isotopy in the family $x^n + y^2 + z^2 = e^{it}$, $t \in [0,2\pi]$, finishes
with the transformation $(x,y,z) \mapsto (\e_n x,-y,-z)$ whose action on the homology coincides with that of $g$. 

\begin{figure}[hbtp]
\begin{center}
\includegraphics[width=0.48\textwidth]{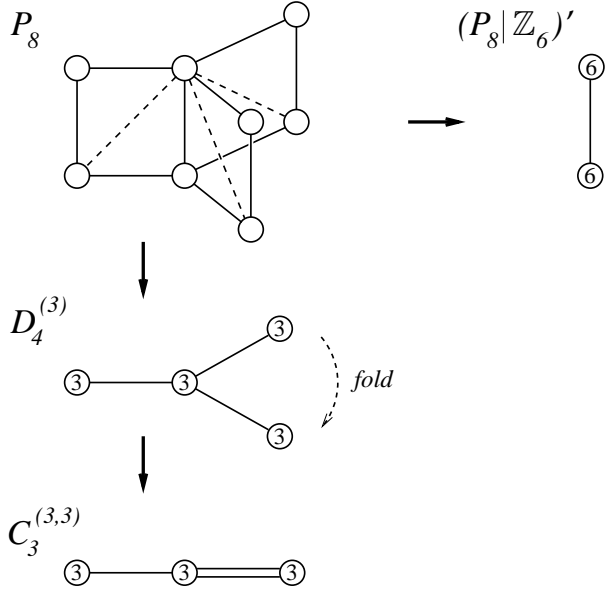}
\caption{Obtaining some of the diagrams by fusion and folding.}
\label{Ffus}
\end{center}
\end{figure}

(b) The $A_2$-version of the previous singularity is $A_2^{(n)}: x^n+y^3+z^2$, with the same symmetry.
For it, each of the $H_{2,\chi}$ is spanned by two similarly defined 
$\chi$-cycles which may be chosen so that their intersection number is $n/(1-\chi)$. The two Picard-Lefschetz operators $h_j$ 
satisfy the standard braiding relation $aba=bab$. 
Diagrammatically, the $A_2^{(n)}$ singularity is represented by fusing the
rectangular $2 \times (n-1)$ Dynkin diagram to the $A_2$ diagram which may be equipped with the markings indicating the orders of the
vertices and the intersection numbers. Cf. the $P_8 \to D_4^{(3)}$ part of Figure \ref{Ffus}. 

(c) The singularity $A_m/\ZZ_m$ is the function $x^{m+1}+yz$ with the equivariant symmetry $g: (x,y,z) \mapsto (\e_m x,\e_m y,z)$. As it has
been shown in \cite{G2}, its $\chi$-cycle has self-intersection $-m$, $\chi$ being any $m$th root of unity. The quasi-homogeneous
argument applied to the $\R_g$-miniversal deformation $x^{m+1}+yz + \a x$ demonstrates that the Picard-Lefschetz operator
on $H_{2,\chi}$ is multiplication by $\chi$.  

\begin{figure}[hbtp]
\begin{center}
\includegraphics[width=0.8\textwidth]{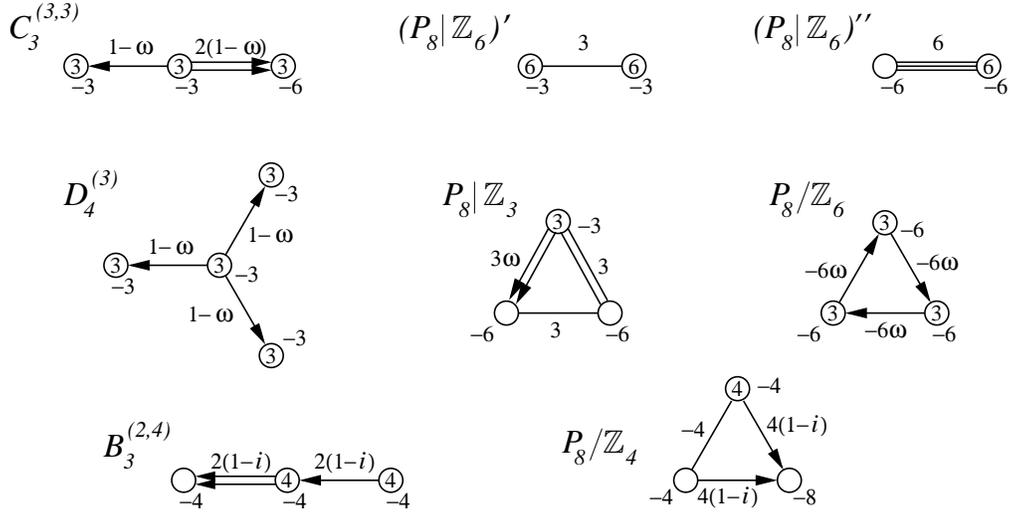}
\caption{Dynkin diagrams of the symmetric $P_8$ singularities. The numbers inside vertices are the orders of the Picard-Lefschetz operators
(order 2 omitted). The number next to a vertex is the self-intersection of the $\chi$-cycle. As earlier, edges are marked with the intersection numbers
of the cycles. In the diagrams of all {\em invariant\/} singularities, the multiplicity of an edge indicates the length of the pair-wise braiding relation
inherited from the fundamental group of the complement to the discriminant:
$aba=bab$ for a simple edge, $(ab)^2=(ba)^2$ for a double and $(ab)^3=(ba)^3$ for a triple.} 
\label{Fdi}
\end{center}
\end{figure}

\begin{thm}\label{Tdd}
For the symmetric $P_8$ singularities and their kernel characters $\chi$, there exist distinguished sets of vanishing 
$\chi$-cycles described by Dynkin diagrams of Figure \ref{Fdi}.
\end{thm}

{\em Proof\/}. We proceed on the case-by-case basis. We are considering only kernel characters of the $P_8$ singularities. 

\medskip
$\mathbf{D_4^{(3)}}$. Factorisation by the action of $\ZZ_3$ provides the boundary $D_4$ singularity. Hence, according to the above
Examples, a Dynkin diagram of $D_4^{(3)}$ can be obtained by the modification of the standard $D_4$ diagram: the 
roots should have squares $-3$ instead of $-2$, the Picard-Lefschetz operators $h_j$ become of order 3, and the intersection
numbers 1 of pairs of cycles change to $3/(1-\chi)$. Since $\chi=\w,\ov\w$, the latter may be reduced to $1-\w$ 
(for both values of $\chi$) using the ambiguities in the choice of the $\chi$-cycles. 
The diagram may be obtained by the fusion of the cylindrical $P_8$ diagram as shown in Figure \ref{Ffus}.   

\medskip
$\mathbf{C_3^{(3,3)}}$. We extend the previous $\ZZ_3$-symmetry by $\ZZ_2$ acting by the sign change on $z$. This embeds the character
subspaces $H_{2,-\w}$ and $H_{2,-\ov\w}$ of the symmetric singularity $C_3^{(3,3)}$ into  
$H_{2,\w}$ and $H_{2,\ov\w}$ of $D_4^{(3)}$ as subspaces anti-invariant under the involution. 
This is absolutely similar to the relation between the $\ZZ_2$-symmetric singularity $C_3$ and the absolute $D_4$ singularity,cf.
\cite{Ab,AGV2,AGLV1}.  
Hence the $C_3^{(3,3)}$ diagram of Figure \ref{Fdi} is provided by the folding of the $D_4^{(3)}$ shown in Figure \ref{Ffus}.

\medskip
$\mathbf{(P_8|\Z_6)'}$. The discriminant in this case is a semi-cubic cusp, hence the relation $aba=bab$ between  
the two Picard-Lefschetz operators. A generic point of the discriminant corresponds to the $D_4|\ZZ_6$ singularity considered in 
\cite{G2}, where the self-intersection number of its vanishing $\chi$-cycle was shown to be $-3$ for $\chi=-\w,-\ov\w$. 
Since we are considering the kernel characters, the intersection form on the $H_{2,\chi}$ must be degenerate which implies that
the absolute value of the intersection number of the two $\chi$-cycles is 3. Since this number is in $\ZZ[1,\w]$, we can make it 3 
using the $\chi$-cycle choice ambiguities once again.
 
Finally, the Picard-Lefschetz operator $h_{D_4|\ZZ_6}$ has order 6 as the classical monodromy operator of the absolute $D_4$ singularity. 
Its non-trivial eigenvalue on $H_{2,\chi}$ is $\l=\ov\chi$ since the quasi-homogeneous isotopy $x^2+y^3+yz^2=e^{it}$, $t\in[0,2\pi]$,
finishes with the transformation $(x,y,z) \mapsto (-x, \w y, \w z)$ which coincides with $g^{-1}$ on the local homology.  

Thus we have arrived at the $(P_8|\ZZ_6)'$ diagram of Figure \ref{Fdi}. Due to the generating reflections coming from the classical 
monodromy of the $D_4$ singularities, it is natural to consider it as obtained by the fusion of Figure \ref{Ffus}.

\medskip
$\mathbf{(P_8|\Z_6)''}$. The discriminant of the versal family $x^3+y^3+yz^2 + \b xy + \a$ consist of two strata:
$$ 
A_5 = \{\a = 0\} \qquad {\rm and} \qquad 3A_1 = \{27\a+\b^3=0\}\,.
$$
The cubic tangency of the strata gives the braiding relation between the generators. 
Of course, the self-intersection of the $3A_1$ vanishing $\chi$-cycle is $-6$ and the order of the reflection $h_{3A_1}$ is 2.
As for the $A_5$ $\chi$-cycle, it is possible to homotope it to the standard $A_1^{(6)}$ vanishing $\chi$-cycle, with the self-intersection
$-6$ according to Example (a). Now the argument similar to that for the previous singularity gives the intersection number $6$ of the
two cycles.  Quasi-homogeneous considerations of the $A_5$ local vanishing form $xy + z^6 = \e$ shows that the non-trivial 
eigenvalue of its monodromy operator on $H_{2,\chi}$ is $\ov\chi$ since its local action coincides with $g^{-1}$. Hence the 
operator $h_{A_5}$ is of order 6. 

\begin{figure}[hbtp]
\begin{center}
\includegraphics[width=0.7\textwidth]{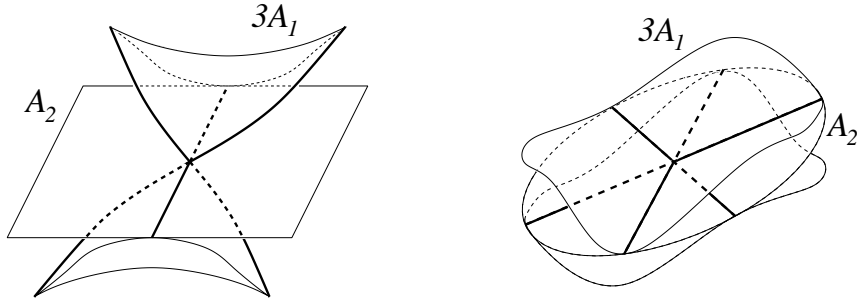}
\caption{Two real versions of the $P_8|\ZZ_3$ discriminant, for $x^3+y^3+yz^2$ and $x^3+y^3-yz^2$.}
\label{Fnew}
\end{center}
\end{figure}

\medskip
$\mathbf{P_8|\Z_3}$. The discriminants of two real versions $x^3+y^3 \pm yz^2 + \a + \b xy + \g xz$ of a versal family are
the unions of the $A_2$ stratum $\a=0$ and $3A_1$ stratum $(54\a+\b^3 \pm 9\b\g^2)^2 = (\b^2 \mp 3\g^2)^3$. They are shown
in Figure \ref{Fnew}. The two strata are simply tangent to each other along their meeting lines.

Like for the previous singularity, it is possible to homotope the $A_2$ vanishing $\chi$-cycle to the standard
$A_1^{(3 )}$ vanishing $\chi$-cycle, with the self-intersection $-3$. An operator $h_{A_2}$ has the eigenvalue $\chi$. 
The self-intersection of a $3A_1$ $\chi$-cycle is $-6$ and its monodromy reflection is an involution. 

Assume the base point $u_*$ is chosen inside the front lips region of the left diagram of Figure \ref{Fnew}. A generic line through it is vertical.
Take on it a path system from $u_*$: two paths going straight to the $3A_1$ points and the third nearly straight to the $A_2$ point, bypassing
the upper $3A_1$ point on its way (the side is not important). Then the vanishing $\chi$-cycles may be chosen so that
$$ 
\left< e_{3A_1,lower},   e_{3A_1,upper} \right> = 3 \qquad {\rm and} \qquad  
\left< e_{A_2},   e_{3A_1,upper} \right> = 3\,.
$$
Indeed the absolute values of these intersection numbers must be 3 following the braiding relations between the
two pairs of the Picard-Lefschetz operators
coming from the singularities of the discriminants. Moreover, both numbers are in $\ZZ[1,\w]$ and the $\chi$-cycles may be multiplied
by $-1$ and powers of $\w$. 


Now we see that the last intersection number, $\left< e_{A_2},   e_{3A_1,lower} \right>$, is either $3\w$ or $3\ov\w$. This is guaranteed
by the degeneracy of the intersection form on $H_{2,\chi}$ and the number being in $\ZZ[1,\w]$. However, both options turn out to be
equivalent up to a braid transformation. 

\medskip
$\mathbf{C_3^{(2,4)}}$. The discriminant coincides with the standard $C_3$ discriminant. Its smooth stratum is $2A_1$
and the cuspidal is $A_3$. The latter is nearly $A_1^{(4)}$ of Example (a), but with the order 4 symmetry changing the sign of one of the squared
variables, hence the self-intersection of a vanishing $\chi$-cycle $e_{A_3}$ is $-4$, while its operator $h_{A_3}$ has eigenvalue $-\ov\chi=\chi$.
Additionally, $\left<e_{2A_1},e_{2A_1}\right>$ is also $-4$ and $h_{2A_1}$ is an involution. 

A distinguished set of vanishing $\chi$-cycles consists of  $e_{2A_1,1},e_{2A_1,2}, e_{2A_1}$. 
The singularities of the discriminant guarantee that the set may be chosen so that the braiding relations between the Picard-Lefschetz
operators are exactly those encoded in the diagram of Figure \ref{Fdi}. These imply 
$$
\left<e_{2A_1,1}, e_{2A_1}\right> = 0 \qquad {\rm and} \qquad 
|\left<e_{2A_1,1},e_{2A_1,2}\right>| = |\left<e_{2A_1,2}, e_{2A_1}\right>| = 2\sqrt2\,. 
$$ 
Since the intersection numbers are in $\ZZ[1,i]$, we can make  
$\left<e_{2A_1,1},e_{2A_1,2}\right> = \left<e_{2A_1,2}, e_{2A_1}\right> = 2(1-i)$ by possible multiplication of the cycles by powers of $i$.

{\em Remark.} We prefer calling this singularity $C$ rather than $B$ since it has only one vanishing cycle coming from more
than one critical point on one level (cf. \cite{Ab}).

\medskip
$\mathbf{P_8/\Z_6}$. The discriminant of the versal family $x^3+y^3+z^3 + \a y + \b z$ consists of 3 lines corresponding to the
three divisors $\a y + \b z$ of the cubic form $y^3 + z^3$. Each of the lines is a $2A_1^{(3)}$ stratum: each $A_1^{(3)}$ singularity
is off the origin and has symmetry $g^2: (x,y,z) \mapsto (\ov\w x,y,z)$. Hence all vanishing $\chi$-cycles have self-intersection $-6$.
Any Picard-Lefschetz operator has eigenvalue $\l = \ov\chi^2 = \chi$.

A distinguished set of $\chi$-cycles consists of 3 elements. 
Since the rank of the intersection form on $H_{2,\chi}$ is 1, the absolute value of the intersection number of any pair is 6.
Since the number itself is in $\ZZ[1,\w]$, we can use the ambiguities in the choice of the cycles and make 
$\left<e_1,e_2\right> = \left<e_2,e_3\right> = -6\w$. The rank condition implies that 
$\left<e_3,e_1\right>$ is also $-6\w$.

{\em Remark.} It is possible to show that the relation between the cycles may be assumed to be $e_1+e_2+e_3=0$. 

\medskip
$\mathbf{P_8/\Z_4}$. The discriminant of the versal family $x^2z + xy^2 + z^3 + \a x + \b yz$ has 3 strata:
$$
A_4: \a = 0\,, \qquad 2A_2: \b = 0 \,, \qquad 4A_1: 4\a + \b^2 =0\,.
$$
The $A_4$ degeneration reduces to the $A_4/\ZZ_4$ singularity of Example (c). Hence the self-intersection of $e_{A_4}$ is $-4$.
By the quasi-homogeneity, we see that the operator $h_{A_4}$ coincides locally with $g$ and thus has eigenvalue $\chi$ on $H_{2,\chi}$.
Similarly, each of the two $A_2$ degenerations is the $g^2$-symmetric singularity $A_2/\ZZ_2$. 
Therefore, the self-intersection of $e_{2A_2}$ is also $-4$, while $h_{2A_2}$ is an involution. It is clear
that  a $\chi$-cycle vanishing at $4A_1$ has the self-intersection $-8$, and $h_{4A_1}$ is also an involution.

A distinguished set of vanishing $\chi$-cycles contains one cycle of each of the three kinds. Since the rank of the $H_{2,\chi}$ is
2, the $3 \times 3$ intersection matrix must have rank 1. Due to this, after possible multiplication of the cycles by powers of $i$,
we can assume that
$$
\left<e_{4A_1},e_{A_4}\right> = \left<e_{4A_1},e_{2A_2}\right> = 4(1+i)\,.
$$
This forces $\left<e_{A_4},e_{2A_2}\right> = -4$. 

{\em Remark.} The relation between the cycles may be assumed to be either 
$e_{A_4}+e_{4A_1}+ ie_{2A_2}=0$ or $ie_{A_4}+e_{4A_1}+ e_{2A_2}=0$.  \hfill{$\Box$}

\section{Affine monodromy}

We start with a reminder how a corank 1 semi-definite hermitian form $\wt q$ on $\wt V = \C^{n+1}$ defines an affine
reflection group on a hyperplane in the dual space \cite{Bou,GM,X9}.

First of all, choose a basis $e_0',e_1,\dots,e_n$ in $\wt V$ so that $e'_0$ is in the kernel $K$ of the form. We denote
the span of the $e_{j>0}$ by $V$ and write $v$ for the $V$-component of $\wt v \in \wt V$: $\wt v = v_0 e_0' + v$. 
Set $Q$ to be the matrix of the restriction $q = \wt q|_V$: $Q=(\wt q (e_i,e_j))_{i,j>0}$.

In the dual space $\wt V^* \simeq K^* \oplus V$, we use coordinates $\a_0,\a_1,\dots,\a_n$ so that a linear functional
$\wt \a$ on $\wt V$ is written as 
$$\wt \a (\wt v) = v_0\a_0 + v^T Q \a = v_0 \a_0 + q(v,\ov \a)\,.$$

Consider a reflection on $\wt V$ with a root $\wt u \notin K$ and the eigenvalue $\l$:
$$
A: \wt v \mapsto \wt v - (1-\l) \wt q (\wt v, \wt u) \wt u / \wt q (\wt u, \wt u)\,.
$$
Then its dual $A^*$ sends each hyperplane $\a_0=const$ in $\wt V^*$ into itself and on such a hyperplane it acts as
$$
\a \mapsto \a - (1 - \ov \l){\a_0u_0 + \ov q(\a,\ov u) \over \ov q(\ov u, \ov u)} \ov u\,,
$$
where $\ov q$ if the hermitian form on $V$ conjugate to $q$, with the matrix $\ov Q = Q^T$ in the chosen coordinates. 
For $\a_0 \ne 0$, this is an affine reflection on the hyperplane $\G = \{\a_0 = const\} \simeq V$, with the root $\ov u$,
mirror $\wt a (\wt u) = \a_0 u_0 + \ov q (\a, \ov u) =0$ and eigenvalue $\ov\l \,.$ For $u_0=0$, the transformation is linear.

\bigskip
We are now ready to prove our main result

\begin{thm} Assume a symmetry of a function singularity $P_8$ on $\C^3$ splits the kernel of the intersection form $\s_h$ 
and $\chi$ is a kernel character. Assume the cohomology character subspace $H^2_{\chi}$ is at least of rank 2.
Let $\G \subset H^2_\chi$ be the hyperplane formed by all 2-cocycles taking a fixed non-zero value on a fixed generator
of the line $H_{2,\chi} \cap Ker(\s_h)$.
Then the equivariant monodromy acts on $\G$ as a complex crystallographic group. The correspondence between the symmetric
singularities and affine groups is given by Table \ref{Ta}.
\end{thm}

{\em Proof\/}.
Number the vanishing $\chi$-cycles of each diagram in Figure \ref{Fdi} from the left to the right
(the numbering of the last two vertices in $D_4^{(3)}$ is not important). 
Omit the leftmost cycle $e_0$. In equivariant cases omit also $e_2$. It is easy to notice that
the remaining cycles may be multiplied by appropriate factors so that the encoded hermitian form becomes the negative 
of the relevant form of Figure \ref{Fgr} while the orders of the related vertices coincide. Therefore,
taking the remaining vanishing cycles  
as the basic vectors $e_{j>0}$ in $\wt V = H_{2,\chi}$ in the
introduction to this Section, we see that their Picard-Lefschetz operators $h_j$
generate the required Shephard-Todd group $L$ on $\G \subset H^2_{\chi} = \wt V^*$. 

The next task is to obtain the translation lattices of the crystallographic groups. 
The kernel of any intersection form in Figure \ref{Fdi} is spanned by the $e_0'=e_0 + {\mathbf a}$, where ${\mathbf a}$ is
a linear combination of the $e_{j>0}$. For example, in all $\tau =2$ cases, $\mathbf a$ is a non-zero multiple of $e_1$
(for equivariant singularities this is due to the two Remarks by the end of the previous Section).


Assume the order of a root of $L$
which generates the lattice of the required affine group coincides with the order of an operator $h_k$, $k>0$. 
According to what has been said before the Theorem, its validity for such singularity will follow from  
$\mathbf a$ being a multiple of an element of the $L$-orbit of $e_k$. Hence only $\tau > 2$ singularities are still to be checked.
And we have for them:
$$
\begin{array}{rcclcl}
C_3^{(3,3)}: & \mathbf a & = & (1-\ov\w, 1) & = & A_1A_2^2 e_1 \\
D_4^{(3)}:    & \mathbf a & = & (1-\ov\w,1,1)             & =  & A_1A_3^2A_1 e_2  \\
P_8|\ZZ_3:    & \mathbf a & = & (\ov\w - \w, -\w)    & =   & \ov\w A_2A_1^2 e_2 \\
C_3^{(2,4)}: & \mathbf a & = & (1+i,i) & = & iA_1^2e_2 
\end{array}
$$
To make the calculations suitable for any of the two kernel characters, the $A_j$ here are either the Picard-Lefschetz
operators or their inverses, but always with the eigenvalues either $\w$ or $-1$ or $i$. The expressions obtained show that
the vectors $\mathbf a$ are maximal roots of the groups $L$ in the sense \cite{H}.
\hfill{$\Box$}

\bigskip
\noindent
{\bf Remark.} The multiplicities of vertices and edges in Figure \ref{Fdi} provide nearly complete presentations of all our rank $>1$ 
crystallographic groups as abstract groups. To obtain all defining relations, one should add:
\begin{itemize}
\item[(i)] $(h_1h_0h_2h_0)^2 = (h_0h_2h_0h_1)^2$ to the $P_8|\ZZ_3$ diagram of $[G(3,1,2)]_2$ (cf. \cite{M}), which corresponds
to one of the tangency lines in Figure \ref{Fnew}, right;
\item[(ii)] the condition that the classical monodromy (that is, the product of all the generators) of each
singularity has order 3. This is exactly the extra relation from \cite{M}.
\end{itemize}

\ 

\begin{tabbing}
Department of Mathematical Sciences,\= \qquad \qquad \qquad\quad\= Department of Mathematics,\\
University of Liverpool, \> \> Ben Gurion University of the Negev,\\
Liverpool L69 7ZL, \> \> P.O.B. 653, Be'er Sheva 84105,\\
United Kingdom\> \> Israel\\
{\em E-mail}: {\tt goryunov@liv.ac.uk} \> \> {\em E-mail}: {\tt kernerdm@math.bgu.ac.il}
\end{tabbing} 

\end{document}